\documentstyle[12pt,fleqn,leqno]{article}
\newtheorem{thm}{Theorem}[section]
\newcommand{\be}{\begin{equation}}
\newcommand{\ee}{\end{equation}}
\newcommand{\ba}{\begin{array}}
\newcommand{\ea}{\end{array}}

\renewcommand{\t}{\theta}

\newcommand{\mg}{\rm}
\renewcommand{\em}{\it}
\newcommand{\bea}{\begin{eqnarray}}
\newcommand{\eea}{\end{eqnarray}}

\newcommand{\Sum}{\sum_{n=1}^\infty}
\newcommand{\Summ}{\sum_{n=0}^\infty}

\begin{document}
\newtheorem{lem}[thm]{Lemma}
\newtheorem{cor}[thm]{Corollary}
\title{\bf A Right Inverse Of The Askey-Wilson Operator \thanks{
Research partially supported  by NSF grant DMS 9203659 and  grants from  SERC.} }
\author{B. Malcolm Brown and Mourad E. H.  Ismail}
\date{}
\maketitle
\begin{abstract}
We establish an integral representation of a right inverse of the Askey-Wilson 
finite difference operator on  $L^2$   with weight $(1-x^2)^{-1/2}$.  The kernel
of this integral operator is $\vartheta'_4/\vartheta_4$ and is the Riemann
mapping function that maps the open unit disc conformally onto the interior
of an ellipse. 
\end{abstract}

\bigskip 
{\bf Running title}:$\;$ An Inverse Operator


\vfill\eject
\section{Introduction.} In 1985 Askey and Wilson  introduced what 
has become known as the 
Askey-Wilson operator. It is defined in the following way. Given a function
 $f(x)$ with $x = \cos\t$ then 
$f(x)$ can be viewed as a function of $e^{i\t}$. Let
\be
\breve{f}(e^{i\t}):=f(x), \quad x = \cos\t.
\ee
In this notation the Askey-Wilson finite difference operator ${\cal D}_q$ is defined by 
\be
({\cal D}_q f)(x) := \frac{(\delta_q \breve{f})(e^{i\t})}{\delta_q \cos\t},
\ee
where
\be
(\delta_q g)(e^{i\t})= g(q^{1/2}e^{i\t}) - g(q^{-1/2}e^{i\t}).
\ee
It follows easily from (1.3) that $\delta_q \cos\t = (q^{1/2} - q^{-1/2})\,i \sin \t$.
 Thus
\be
({\cal D}_q f)(x) = \frac{\breve{f}(q^{1/2}e^{i\t}) - \breve{f}(q^{-1/2}e^{i\t})}
{i(q^{1/2}- q^{-1/2})\sin \t}.
\ee

The operator ${\cal D}_q$ was introduced in \cite{As:Wi} and in the theory of the
Askey-Wilson polynomials plays an analogous role  to that of differentiation in 
the theory of Jacobi, Hermite and Laguerre polynomials.  
 Alphonse Magnus \cite{Ma} provided a  general setting for 
finite difference operators and indicated how the Askey-Wilson
 operator arises in a very natural way as a finite difference operator. 
We remark that 
${\cal D}_q$ remains invariant if $q$ is replaced by $1/q$. In this work we
will avoid $q$'s on the unit circle, thus there is no loss of generality in
assuming $|q| < 1$.

 The Chebyshev polynomials $T_n$ and $U_n$ of the first and second kinds,
 respectively, are
\be
T_n(x) = \cos n\t, \quad U_n(x) = \frac{\sin ((n+1)\t)}{\sin \t}, 
\quad x = \cos\t,
\ee
 and are orthogonal on $[-1, 1]$ with respect to
the weight functions $(1 - x^2)^{-1/2}$ and $(1 - x^2)^{1/2}$, respectively. A 
calculation gives
\be
{\cal D}_q T_n(x) = \frac{q^{1/2} - q^{-n/2}}{q^{1/2} - q^{-1/2}}\, U_{n-1}(x).
\ee
The Askey-Wilson polynomials $p_n(x; a, b, c, d|q)$ \cite{As:Wi}, \cite{Ga:Ra}, 
which depend on four parameters, are orthogonal 
on $[-1, 1]$ when all their parameter values are in $(-1, 1)$. Thus one would 
like to study the action of ${\cal D}_q$ on
 functions defined on $[-1,1]$. One difficulty in defining the action of
 ${\cal D}_q$,  defined in (1.4), is that (1.4) uses values of $f$ outside 
$[-1, 1]$.
 In fact the definition (1.4) given in \cite{As:Wi} uses values of $f$ at  points
 in the complex plane
  without specifying how $f$ is extended to the  whole, or part of,
the $x$-plane. To rectify these difficulties we propose to define 
${\cal D}_q$ on a dense subset of $L^2[(1-x^2)^{-1/2}, [-1, 1]]$. We propose
 the following definition.

{\bf Definition 1.1}. Let $f \in L^2[(1-x^2)^{-1/2}], [-1, 1]]$, and assume that
$f$ has a Fourier-Chebyshev expansion
\be
f(x) \sim \Summ f_n T_n(x),
\ee
with the Fourier-Chebyshev coefficients  $\{f_n\}$ satisfying
\be
\Summ |(1-q^n)\, q^{-n/2}\, f_n|^2 < \infty.
\ee
We then define ${\cal D}_q f$ as the unique (almost everywhere) function 
whose Fourier-Chebyshev expansion is
\be
({\cal D}_q f)(x) \sim \Summ \frac{q^{n/2} - q^{-n/2}}{q^{1/2} - q^{-1/2}}
\, f_n\,U_{n-1}(x).
\ee
Such functions $f$ we call $q$-differentiable.

It is clear that (1.7) is satisfied on a dense subset $S$ of 
$L^2[(1-x^2)^{-1/2}, [-1, 1]]$. It is also clear that ${\cal D}_q$ maps $S$ 
into $L^2[(1-x^2)^{1/2}, [-1, 1]]$.

In Section 2 we define a right inverse to ${\cal D}_q$, which will be denoted as
${\cal D}_q^{-1}$, and find an integral
representation for this inverse operator.  We shall show that the 
 operator ${\cal D}_q^{-1}$ is 
a convolution operator whose kernel is $\vartheta'_4/\vartheta_4$. Also in Section 2 we prove that  ${\cal D}_q
{\cal D}_q^{-1}$ is the identity operator on $L^2[(1-x^2)^{-1/2}, [-1, 1]]$.
Further we  prove that if $g(x)$ is a piecewise continuous  function on $[-1, 1]$ 
with ${\cal D}_q$ defined by (1.4), then $({\cal D}_p{\cal D}_q^{-1}g)(x)$ converges to $g(x)$ as $p \to q^+$ at 
the points of continuity of $g$. Section 3 contains remarks on related work and
establishes a connection between the kernel of ${\cal D}_q^{-1}$ and the
 Riemann mapping 
function of the interior of an ellipse to the open unit disc. The definition 
(1.9) of ${\cal D}_q$ uses $\{T_n(x)\}$ as a basis for a function space. 
In Section 4 we use a different basis, namely the  Rogers $q$-Hermite polynomials 
in the $L^2$ spaces weighted by their weight function. We again identify
 ${\cal D}_q^{-1}$ as an integral operator and find its kernel explicitly.

\section{An Integral Operator.}
\setcounter{equation}{0}
We seek an operator ${\cal D}_q^{-1}$ so that
\be
{\cal D}_q {\cal D}_q^{-1} = I.
\ee
Let ${\cal D}_q f = g$ so that $f(x) \sim \Summ f_n\,T_n(x)$,
$g(x) \sim \Sum g_n\, U_{n-1}(x)$ and
\be
f_n = g_n\, \frac{(q^{1/2} - q^{-1/2})}{(q^{n/2}-q^{-n/2})}, \quad n > 0,
\ee
we first use a heuristic (or formal) approach to find a way to recover $f$
from the knowledge of $g$. Formally we have 
\bea
\Sum f_n\, T_n(x)  &=& \Sum g_n \, \, \frac{(q^{1/2} - 
q^{-1/2})}{(q^{n/2}-q^{-n/2})} T_n(x)
\nonumber \\
& = & \frac{2}{\pi}\, 
\Sum \left(\int_{-1}^{1}g(y)U_{n-1}(y)\,\sqrt{1-y^2}\, dy\right)
\, \frac{(q^{1/2} - q^{-1/2})}{(q^{n/2}-q^{-n/2})}\, T_n(x)
\nonumber \\
& = &\frac{2}{\pi} (1-q)q^{-1/2}\, \int_{-1}^{1} g(y) \left[\Sum
\frac{T_n(x)U_{n-1}(y)}{1 -q^n}q^{n/2}\right]\sqrt{1 - y^2}dy. \nonumber 
\eea
This suggests defining ${\cal D}_q^{-1}$ as an integral operator whose kernel 
is
\be
F(x,y) := \frac{2(1-q)}{\pi\sqrt{q}}\Sum \frac{T_n(x)U_{n-1}(y)}{1 -q^n}q^{n/2}.
\ee
Set
\be
x = \cos\t, \quad y = \cos\phi,
\ee
to get
\be
F(\cos\t, \cos\phi) = \frac{(1-q)q^{-1/2}}{\pi\sin\phi} \Sum
\frac{2\cos(n\t)\, \sin(n\phi)}{1 - q^n}q^{n/2}
\ee
\bea
=  \frac{(1-q)q^{-1/2}}{\pi\sin\phi} \Sum
\frac{q^{n/2}}{1 - q^n} [\sin (n(\t + \phi)) - \sin (n(\t-\phi))],
\nonumber
\eea
and observe that 
\bea
\int_{-1}^1 F(x, y)g(y)\sqrt{1-y^2} dy &  =  & \int_0^{\pi} F(\cos\t, \cos\phi)
g(\cos\phi)\sin^2\phi d\phi  \nonumber 
\\ & = & \int_{-\pi}^{\pi}G(\cos\t, \cos\phi)
 g(\cos\phi)\sin^2\phi d\phi,
\nonumber
\eea
where
\be
G(\cos\t, \cos\phi) = \frac{(1-q)}{\pi\sqrt{q}}\Sum \frac{q^{n/2}}{1-q^n}
\sin (n(\t+\phi)).
\ee
We next recall that the theta function $\vartheta_4(z, q)$ may be defined by
\be
\vartheta_4(z, q) = \sum_{-\infty}^{\infty} (-1)^nq^{n^2}e^{2inz},
\ee
see \cite[p. 463]{Wh:Wa}, the logarithmic derivative of which  has the Fourier
 series expansion
\be
\frac{\vartheta'_4(z, q)}{\vartheta_4(z, q)} = 4 \Sum \frac{q^n}{1-q^{2n}}
\sin (2nz),
\ee
 \cite[Ex.11, p. 489]{Wh:Wa}. Thus (2.3), (2.4), (2.5) and (2.6) motivate the
 following definition:

\noindent{\bf Definition 2.1}. The operator ${\cal D}_q^{-1}$ is defined on 
 $L^2[(1 - x^2)^{1/2}, [-1,1]]$ as the integral operator
\be
({\cal D}_q^{-1}g)(\cos \t) = \frac{1-q}{4\pi \sqrt{q}} \int_{-\pi}^{\pi}
\frac{\vartheta'_4((\t-\phi)/2, \sqrt{q})}{\vartheta_4((\t-\phi)/2, \sqrt{q})}
\, g(\cos \phi) \, \sin \phi\, d\phi.
\ee

Observe that the kernel of the integral operator (2.9) is bounded when $(x, y)
 (=(\cos \t, \cos \phi)) \in [-1, 1] \times [-1,1]$. Thus the operator
 ${\cal D}_q^{-1}$ is well-defined and 
bounded on $L^2[(1-x^2)^{1/2}, [-1,1]]$. Furthermore ${\cal D}_q^{-1}$  is a
one-to-one mapping from
 $L^2[(1-x^2)^{1/2}, [-1,1]]$ into  $L^2[(1-x^2)^{-1/2}, [-1,1]]$. 

\begin{thm} The operator  ${\cal D}_q {\cal D}_q^{-1}$ is 
the identity operator on $L^2[(1-x^2)^{1/2}, [-1,1]]$. 
\end{thm}

{\bf Proof}. Replace $\vartheta'_4/\vartheta_4$ in (2.9) by the expansion (2.8)
 then apply Parseval's formula. In view of the uniform convergence of the series
in (2.8) we may  reverse the order of integration and summation in (2.9) 
from which it follows that the  steps leading to  the Definition
2.1  can now be reversed  yielding the result.

Note that the kernel $G$ of (2.6) is defined and bounded for all $\phi \in [-\pi, \pi]$ and  all $\t$ for which $|q^{-1/2} e^{i\t}| \le 1$, where we have written 
 $e^{i\t} := x +
 \sqrt{x^2-1}$ and the branch of the square root  is chosen so that 
 $\sqrt{x^2-1} \approx x$ as $x \to \infty$. This convention makes 
$|e^{-i\t}| \le |e^{i\t}|$
with equality if and only if $\t$ is real, that is $x \in [-1, 1]$.  
This extends the definition of ${\cal D}_q^{-1}$ to the interior of the
 ellipse $|z + \sqrt{z^2 - 1}| = q^{-1/2}$ in the complex $z$-plane. This 
ellipse has foci at $\pm1$ and its major and minor axes are 
$q^{-1/2} \pm q^{1/2}$, respectively. Its equation in the $xy$-plane is
\be
\frac{x^2}{a^2} + \frac{y^2}{b^2} = 1, \quad a =  \frac{1}{2}(q^{-1/2}+ q^{1/2}), 
\quad  b = \frac{1}{2}(q^{-1/2} - q^{1/2}).
\ee

It is worth mentioning that provided we exercise some care the Askey-Wilson
 definition (1.4) yields the result
${\cal D}_q {\cal D}_q^{-1} = I$. One
reason  for being particularly carefulll is that ${\cal D}_q^{-1}g$ may not be
 in the domain of 
${\cal D}_q$
because in order to use (1.3) we need to assume that $f$ has an analytic
extension to a domain  in the complex plane containing $|z \pm \sqrt{z^2-1}| 
\le q^{-1/2}$.  Let $f(\cos \t)$ denote the  right-hand side of (2.9).  From the
 discussion following Theorem 2.1 it is easy to see that ${\cal D}_p f(x)$ is well-defined provided that $1 > p > q$. Indeed
we find
\bea
({\cal D}_p f) (\cos \t) &=& \int_{-\pi}^{\pi}{\cal D}_p G(\cos \t, \cos \phi)\, \sin^2 \phi \, d\phi  \nonumber \\
& = & \frac{(1-q)q^{-1/2}}{\pi \sin \t} \int_{-\pi}^{\pi} \Sum  \frac{(1-p^n)}
{(1-q^n)}\,(q/p)^{n/2} \cos (n(\t - \phi))\, g(cos \phi)\, \sin \phi \, d\phi. \nonumber 
\eea
By writing $\frac{1-p^n}{1-q^n}$ as $1 + \frac{q^n - p^n}{1 - q^n}$, and denoting
$q/p$  by $r$ we see that
\bea
\lim_{p \to q^+} ({\cal D}_p{\cal D}_q^{-1}g)(\cos t) &=& \frac{1}
{\pi \sin \t} \lim_{r \to 1^-} \int_{-\pi}^{\pi}\left[\Sum r^n
\cos (n(\t - \phi))\right]\, g(\cos \phi)\, \sin \phi \, d\phi  + 0 \nonumber \\
&=& \frac{(1-q)q^{-1/2}}
{\pi \sin \t} \lim_{r \to 1^-} \int_{-\pi}^{\pi}\left[ \frac{1}{2} + \Sum r^n
\cos (n(\t - \phi))\right]\, g(\cos \phi)\, \sin \phi \, d\phi, \nonumber 
\eea
since $g(\cos \phi) \sin \phi$ is an odd function. Therefore
\bea
\lim_{p \to q^+} ({\cal D}_p{\cal D}_q^{-1}g)(\cos t) = 
\lim_{r \to 1^-} \frac{1}{2\pi \sin \t} \int_{-\pi}^{\pi}\frac{(1-r^2)
\, g(\cos \phi)\, \sin \phi}{1 - 2r\cos (\t - \phi) +r^2}\, d\phi. \nonumber
\eea
Finally the above limit exits and equals $g(\cos \t)$ at the points of
 continuity  of $g$ if $g(\cos \t)$ is continuous on $[-\pi, \pi]$ except 
for finitely many jumps, \cite[p. 147]{Ne2}. Thus we have proved the following 
result.
\begin{thm}
Let $g(x)$ be a continuous function on $[-1,1]$ except for finitly many jumps.
Then with ${\cal D}_q$ defined as in (1.4), the limiting relation
\bea
\lim_{p \to q^+} ({\cal D}_p{\cal D}_q^{-1}g)(x) = g(x)  \nonumber
\eea
holds at the points of continuity of $g$.
\end{thm}

\section{Remarks.} 
\setcounter{equation}{0}
 The kernel $\vartheta'_4/\vartheta_4$ of (2.8) has appeared earlier in conformal
 mappings. Let $\zeta$ be a fixed point in the interior of 
the ellipse (2.10)  in the complex plane and let $f(z, \zeta)$ be the function
that maps the interior of the ellipse (2.10) conformally onto the open unit disc
 and satisfies $f(\zeta, \zeta) = 0$ and $f'(\zeta, \zeta) > 0$. It is known,
\cite[p. 260]{Ne1}, that
\be
f(z, \zeta) = g(z, \zeta) - g(\zeta, \zeta),
\ee
where
\be
g(z, \zeta) = \sqrt{\frac{\pi}{K(\zeta, \zeta)}} \Sum \frac{T_n(z)\, 
\overline{U_n(\zeta)}}{\rho^n - \rho^{-n}},
\ee
and the Bergman kernel $K(z, \zeta)$ is
\be
K(z, \zeta) := \frac{4}{\pi} \Summ \frac{(n+1)\, U_n(z) \, 
\overline{U_n(\zeta)}}
{\rho^{n+1} - \rho^{-n-1}}, \quad \rho:= (a + b)^2 = (b + \sqrt{b^2+1})^2.
\ee
In fact the Bergman kernel of the ellipse is a constant multiple of 
$f'(z, \zeta)$. It is clear that $g(z, \zeta)$ is a constant multiple of our
 kernel $G$ (2.6) with $\rho = q^{-1/2}$, so that $q = (b + \sqrt{b^2+1)})^{-4} 
= e^{-4u}$ if
$ b = \sinh u$.

The connection between the Riemann mapping function $f(z, \zeta)$ of the 
ellipse (2.10) and our kernel may seem very surprising at a first glance. 
However  this may not be a complete surprise because if $f$ 
is real analytic in $(-1, 1)$ it will have an extension which is analytic 
in the open unit disc and (1.4) will be meaningful if $|q^{-1/2} e^{i\t}| < 1$;
 which is the interior of the ellipse (2.10). Furthermore the Chebyshev
polynomials $\{U_n(z)\}$ are orthogonal on the unit disc with respect to the 
Lebesgue measure in the plane.

Ismail and Zhang  \cite{Is:Zh} proved that the eigenvalues of the integral 
operator (2.9)  are $\pm i/j_{0,k}(q)$ where $\{j_{\nu,k}(q)\}$ 
are the zeros of the $q$-Bessel function $J_{\nu}^{(2)}(z;q)$, see \cite{Is} or
 \cite{Ga:Ra} for the definition of $J_{\nu}^{(2)}(z;q)$. Ismail and Zhang
also proved that the eigenfunctions of ${\cal D}_q^{-1}$ provide a new
 $q$-exponential function. In addition they found the eigenvalues and
 eigenfunctions of ${\cal D}_q^{-1}$ on $L^2$ spaces weighted by the the weight
functions of the Jacobi and continuous $q$-ultraspherical polynomials, \cite{As:Is}, \cite{As:Wi}. Later
Ismail, Rahman and Zhang \cite {Is:Ra:Zh} extended this investigation to the
weights of continuous $q$-Jacobi polynomials.

\section{The $q$-Hermite Space}
\setcounter{equation}{0}
It is a consequence of our definition of ${\cal D}_q$ that that subset $S$
of $L^2[(1-x^2)^{-1/2}], [-1,1]]$ of $q$-differentiable functions gets mapped into  $L^2[(1-x^2)^{1/2}], [-1,1]]$.  In this section we shall define 
 ${\cal D}_q$ on a different weigted $L^2$ space in such a way that 
$q$-differentiable functions are mapped into the same space.

Recall that the definition of the $q$-shifted factorial, \cite{Ga:Ra}
\be
(a;\;q)_0:=1,\quad (a;\;q)_n:=\prod_{j=1}^{n}(1-aq^{j}),\quad n=1,2,\ldots,
\mbox{ or }\infty,
\ee
and the multiparameter notation
\be
(a_1,\ldots,a_m;\:q)_n=\prod_{k=1}^{m}\,(a_k;\;q)_n.
\ee
The $q$-Hermite polynomials of L. J. Rogers are, \cite{As:Is}, \cite{Ga:Ra},
\be
H_n(x|q) := \sum_{k = 0}^n \frac{(q;q)_n}{(q;q)_k(q;q)_{n-k}}\, e^{i(n-2k)\t},
\quad x = \cos \t.
\ee
The orthogonality relation of the $H_n$'s is, \cite{As:Is}, \cite{Ga:Ra},
\be
\int_{-1}^1H_m(x|q) H_n(x|q) w(x) dx = \frac{2\pi (q;q)_n}{(q; q)_\infty}
\, \delta_{m,n},
\ee
where the weight function $w$ is
\be
w(x) = (1 - x^2)^{-1/2} \prod_{n=0}^{\infty}(1 - 2(2x^2-1)q^n + q^{2n}),
\ee
that is $w(\cos \t) = (e^{2i\t}, e^{-2i\t}; q)_\infty /\sin \t$.
The $H_n$'s have the generating function
\be
\Summ \frac{H_n(x|q)}{(q; q)_n}r^n = \; 1/(re^{i\t}, re^{-i\t}; q)_\infty.
\ee
By applying ${\cal D}_q$ to both sides of (4.6) we get
\be
{\cal D}_q H_n(x|q) = 2q^{(1-n)/2} \frac{(1-q^n)}{(1-q)}\; H_{n-1}(x|q).
\ee
One can then follow the procedure of Section 2 except that we now work with
 the function space $L^2[w(x), [-1, 1]]$ and ${\cal D}_q$ maps it into itself.
We then establish, in a straight forward manner, the integral representation
\be
({\cal D}_q^{-1}g)(x) =\frac{(1-q)}{4\pi \sqrt{q}}\int_{-1}^1 H(x, t, \sqrt{q})
 \,g(t)\, w(t)\, dt,
\ee
where the more general kernel $H(x, t, r)$  is defined to be
\be
H(x, t, r) = \Sum \frac{r^r}{(q;q)_n}\;H_n(x|q)\, H_{n-1}(x|q).
\ee
We now evaluate the kernel $H(x, t, r)$.
\begin{thm}
We have
\bea
H(\cos \t, \cos \phi, r) &=& \frac{r\,(qr^2; q)_\infty}{(re^{i(\t+\phi)}, 
re^{i(\t-\phi)}, re^{-i(\t+\phi)}, re^{-i(\t-\phi)};q)_\infty}\\
&\mbox{}& \times \Summ \frac{(re^{i(\t+\phi)}, 
re^{i(\t-\phi)}, re^{-i(\t+\phi)}, re^{-i(\t-\phi)};q)_k}{(qr^2;q)_k}
\left[-cos \t + rq^k \cos \phi\right]. \nonumber
\eea
\end{thm}
 {\bf Proof}. The Poisson kernel of the Rogers $q$-Hermite polynomials is
\be
\Summ \frac{H_n(\cos \t|q) H_n(\cos \phi|q)}{(q;q)_n}r^n 
\ee
\bea
\quad =
\frac{(r^2; q)_\infty}{(re^{i(\t+\phi)}, 
re^{i(\t-\phi)}, re^{-i(\t+\phi)}, re^{-i(\t-\phi)};q)_\infty}.\nonumber
\eea
 Keeping $\t$ fixed and applying ${\cal D}_q$, acting on the variable
 $\cos \phi$, to (4.12) gives
\be
\Sum \frac{H_n(\cos \t|q) H_{n-1}(\cos \phi|q)}{(q;q)_n}r^n q^{-n/2}
\ee
\bea
= \frac{(r^2; q)_\infty\; rq^{-1/2}(-\cos \t + rq^{-1/2}\cos \phi)}{(rq^{-1/2}e^{i(\t+\phi)}, 
rq^{-1/2}e^{i(\t-\phi)}, rq^{-1/2}e^{-i(\t+\phi)}, rq^{-1/2}e^{-i(\t-\phi)};q)_\infty}.\nonumber
\eea
This implies
\bea
H(\cos \t, \cos \phi, rq^{-1/2}) - H(\cos \t, \cos \phi, rq^{1/2}) \nonumber
\eea
\bea
= 
\frac{(r^2; q)_\infty\; rq^{-1/2}(-\cos \t + rq^{-1/2}\cos \phi)}{(rq^{-1/2}e^{i(\t+\phi)}, 
rq^{-1/2}e^{i(\t-\phi)}, rq^{-1/2}e^{-i(\t+\phi)}, rq^{-1/2}e^{-i(\t-\phi)};q)_\infty}.  \nonumber
\eea
Next replace $r$ by $rq^{k+1/2}$ and add the resulting equations. In view of 
$H(x, y, 0) = 0$ the result is
\bea
H(x, y, r) = \frac{r\,(r^2q^{2k+1}; q)_\infty}{(rq^ke^{i(\t+\phi)}, 
rq^ke^{i(\t-\phi)}, rq^ke^{-i(\t+\phi)}, rq^ke^{-i(\t-\phi)};q)_\infty}. \nonumber
\eea
The above equation is equivalent to (4.10) and the proof is complete.

Note that the right-hand side of (4.10) can be written as a sum of two
$_5\phi_4$ functions.

{\bf Acknowledgements}. Thanks to Richard Askey for comments and encouragement
and to Ruiming Zhang for reference \cite{Ne1}.

Department of Computing Mathematics, University of Wales College of Cardiff, 
Mathematics Institute, Senghennydd Road, Cardiff CF2 4YN, United Kingdom

Department of  Mathematics, University of South Florida, Tampa, Florida, 
33620, USA.
\end{document}